\documentstyle[amstex, 11pt]{amsart}
\makeatletter
\title[ideals of minimal mixed multiplicity]{On the depth of blow-up rings of ideals of minimal mixed multiplicity}
\author{ Clare D'Cruz} 
\address{Chennai Mathematical Institute, 92, G.N. Chetty Road, Chennai 600 017,
India.}
\email{clare@cmi.ac.in}

\newcommand{\ncom}{\newcommand}
\ncom{\bq}{\begin{equation}}
\ncom{\eq}{\end{equation}}
\ncom{\beqn}{\begin{eqnarray*}}
\ncom{\eeqn}{\end{eqnarray*}}
\ncom{\beq}{\begin{eqnarray}}
\ncom{\eeq}{\end{eqnarray}}
\ncom{\been}{\begin{enumerate}}
\ncom{\eeen}{\end{enumerate}}
\ncom{\nno}{\nonumber}
\ncom{\hs}{\mbox{\hspace{.25cm}}}
\ncom{\rar}{\rightarrow}
\ncom{\lrar}{\longrightarrow}
\ncom{\Rar}{\Rightarrow}
\ncom{\noin}{\noindent}
\newtheorem{thm}{Theorem}[section]
\newtheorem{lemma}[thm]{Lemma}
\newtheorem{cor}[thm]{Corollary}
\newtheorem{pro}[thm]{Proposition}
\newtheorem{example}[thm]{Example}
\newtheorem{remark}[thm]{Remark}
\newtheorem{definition}[thm]{Definition}
\ncom{\bt}{\begin{thm}}
\ncom{\et}{\end{thm}}
\ncom{\bl}{\begin{lemma}}
\ncom{\el}{\end{lemma}}
\ncom{\bco}{\begin{cor}}
\ncom{\eco}{\end{cor}}
\ncom{\bp}{\begin{pro}}
\ncom{\ep}{\end{pro}}
\ncom{\bex}{\begin{example}}
\ncom{\eex}{\end{example}}
\ncom{\brm}{\begin{remark}}
\ncom{\erm}{\end{remark}}
\ncom{\bdef}{\begin{definition}}
\ncom{\edeff}{\end{definition}}
\ncom{\bc}{\begin{center}}
\ncom{\ec}{\end{center}}
\ncom{\comx}{I\!\!\!\!C}
\ncom{\zee}{$Z\!\!\!\!Z$}
\ncom{\ze}{Z\!\!\!\!Z}
\ncom{\Q}{$I\!\!\!\!Q$}
\ncom{\N}{I\!\!N}
\ncom{\sz}{\scriptsize}
\ncom{\sop}{system of parameters}
\ncom{\eop}{\hfill{$\Box$}}
\ncom{\tfae}{the following are equivalent:}
\ncom{\mm}{minimal multiplicity }
\ncom{\f}{\frac}
\ncom{\la}{\lambda}
\ncom{\si}{\sigma}
\ncom{\ssize}{\scriptsize}
\ncom{\al}{\alpha}
\ncom{\be}{\beta}
\ncom{\Si}{\Sigma}
\ncom{\ga}{\gamma}
\ncom{\kbar}{\overline{\kappa}}
\ncom{\bib}{\bibitem}
\ncom{\sst}{\subset}
\ncom{\sms}{\setminus}
\ncom{\seq}{\subseteq}
\ncom{\est}{\emptyset}
\ncom{\bighs}{\hspace{.5 cm}}
\ncom{\ulin}{\underline}
\ncom{\olin}{\overline}
\ncom{\bip}{\bigoplus}
\ncom{\sta}{\stackrel}

\def\m{{\frak m}}

\def\a{{\bf a}}
\def\depth{\mathop{\rm depth}}
\def\F{{\mathcal F}}

\newcommand{\mptm}{\fontfamily{ptm}\fontsize{10}{12}\selectfont}
\newcommand{\mptmcomx}{\fontfamily{ptm}\fontsize{9}{11}\selectfont}

\begin{document} 

\maketitle
\mptm

\begin{abstract}
We show that if $(R, \m)$ 
is a Cohen-Macaulay local ring and 
$I$ is an ideal of minimal mixed multiplicity, then $\depth~G(I) \geq d- 1$ implies that 
$\depth~F(I) \geq d-1$. We use this to show that if $I$ is a contracted ideal in a two dimensional regular local ring then 
$\depth~R[It]-1= \depth~G(I) = \depth~F(I)$. We also give an infinite class of ideals where $R[It]$ is Cohen-Macaulay but $F(I)$ is not.
\end{abstract}

\section{Introduction}
Throughout this paper we will assume that $(R,\m )$ is a local ring of
dimension $d$ with infinite residue field.  
Contracted ideals in a two dimensional regular local ring have been studied by Zariski \cite{zar}. 
A special class of contracted ideals are complete ideals. Zariski's work reveals that complete ideals 
satisfy several interesting properties. We would like to know if some of these results hold good for contracted ideals.

Associated to any  ideal $I$ in  local ring $(R, \m)$ we have the blowup rings,
namely the the Rees ring $R(I) :=R[It] = \oplus_{n \geq 0} I^n t^n$, the associated graded ring 
$G(I) := R(I) \otimes_R R/I = \oplus_{n  \geq 0} I^n/I^{n+1}$  and the fiber cone 
$F(I) := \oplus_{n   \geq 0} I^n/ \m I^n$. 
It is well known that all the blowup rings associated to a complete ideal in a 
two dimensional regular local ring are Cohen-Macaulay. 
One would naturally address the same question for contracted ideals. 
An interesting result which plays an important role in the Cohen-Macaulayness of blowup rings 
is the well known result of Goto and Shimoda which gives a necessary and sufficient condition 
for the Rees ring to be Cohen-Macaulay \cite{goto-s}. From their result it follows  that if 
$(R, \m)$ is a Cohen-Macaulay local ring, $I$ is an $\m$-primary ideal, and the Rees ring is Cohen-Macaulay, 
then so is the associated graded ring.  
By using the Briancon-Skoda  theorem, it was shown
 that for a  regular local ring $(R,\m)$, the ring  $R(I)$ is
 Cohen-Macaulay  if and only if $G(I)$ is Cohen-Macaulay for any
 $\m$-primary ideal \cite{huneke}. We have extended this result of Huneke to relate 
 the depth of the Rees ring and the associated graded ring of and $\m$-primary ideal in 
 a regular local ring of positive dimension (see Theorem~\ref{rees-graded}).

Motivated by these results,  a natural question to ask is whether the depth of the fiber cone is 
related with the depth of the other blowup rings. For complete ideals in a regular local ring, 
all the blowup rings are Cohen-Macaulay. In \cite{jayanthan}, it was shown that if 
$I$ is a lex-segment ideal in a two-dimensional polynomial ring, then $\depth(R(I)) = \depth(G(I)) = \depth(F(I))$.
In this paper we show that this result holds true  for  any $\m$-primary contracted ideal in a two dimensional regular local ring (see Theorem~\ref{rees-graded-fiber}).

Contracted ideals in a two dimensional regular local ring are precisely ideals of minimal mixed multiplicity (Lemma~\ref{cont-mmm}). For   an $\m$-primary contracted ideal $I$  in a two dimensional regular ring $(R, \m)$, 
$\depth(R(I)) -1= \depth(G(I)) = \depth(F(I))$ (see Theorem~\ref{rees-graded-fiber}). Since $\m$-primary contracted ideals in a two dimensional regular local ring are ideals of minimal mixed multiplicity, we study ideals of minimal mixed multiplicity in a Cohen-Macaulay local ring of dimension at least two. For ideals of minimal mixed multiplicity we show that 
 if $G(I)$ is Cohen-Macaulay, then $\depth(F(I)) \geq d-1$ (see Proposition~\ref{depth}). 
 
 Several examples of $G(I)$ being Cohen-Macaulay but $F(I)$ is not is known, but an example of the Rees ring being Cohen-Macaulay but $F(I)$ is not was not known. In Example~\ref{ex-1} we produce an infinite class of complete ideals in a regular local ring of dimension three where the ideal  Rees ring is Cohen-Macaulay but fiber cone is not.
 The next question of interest is: If $F(I)$ is Cohen-Macaulay, is $G(I)$ and $R(I)$ Cohen-Macaulay. This question has a positive answer when $I$ is an in ideal of minimal mixed multiplicity (see Theorem~\ref{fiber-cm}).

I am very  grateful to J. K. Verma who motivated me into this question.
 I also thank I.~Swanson for encouraging me to write this
paper.

%during a conference  which rekindled my interest. I
%also thank W. Heinzer for reading the original manuscript and for his
%suggestions.  

\section{Depth  of $R(I)$, $G(I)$ and $F(I)$}

In this section we first prove Theorem~\ref{rees-graded} which  relates the depth of the Rees ring and the associated graded ring of $\m$-primary ideals in a regular local ring. Our main result is for contracted ideals. We first show for and $\m$-primary ideal in a two dimensional regular local ring, Cohen-Macaulayness of $R(I)$ implies the Cohen-Macaulayness of $G(I)$ and $F(I)$. We show that if $I$ is an $\m$-primary ideal in a two dimensional regular local ring then $G(I)$ is Cohen-Macaulay if and only if 
$\depth(F(I))$ is Cohen-Macaulay.

%A nice class of ideals where the converse of Theorem~\ref{rees-fiber} holds true is complete ideals. 

Let $I$ be an ideal in a regular local ring $(R, \m)$. The integral closure of an ideal, denoted by
${\olin I}$ is defined as 
${\olin I} 
= \{x \in R| x^n + a_1 x^{n-1} + \cdots + a_n=0; a_i \in I^i\}$.
An ideal $I$ is complete  if $I={\olin I}$. 

\bt
Let $I$ be a complete ideal in a two dimensional regular local ring. Then 
$R(I)$, $G(I)$ and $F(I)$ are all Cohen-Macaulay.
\et
\pf A complete ideal has $r(I)=1$. It is well known that $G(I)$ is Cohen-Macaulay. Hence   $R(I)$ is also Cohen-Macaulay \cite{goto-s}. Cohen-Macaulayness of 
$F(I)$ follows from \cite{shah}. 
\qed

We now study the depth of contracted ideals in a two-dimensional regular local ring. 
An ideal $I$ is contracted if there exists $S= R[m/x]$ such that 
$I = IS \cap R$. In a two-dimensional regular local ring an
$m$-primary complete ideal is always contracted \cite{zar}. We
prove this result 
for any dimension $d$.

We begin with a more general result which relates the depth of the Rees ring and the associated graded ring.

\bt
\label{rees-graded}
Let $I$ be an $\m$-primary ideal in a regular local ring. Then 
$\depth(R(I)) = \depth(G(I))+1$.
\et
\pf
Since $R$ is a  regular local ring  the ring  $R(I)$ is
 Cohen-Macaulay  if and only if $G(I)$ is Cohen-Macaulay for any
 $\m$-primary ideal \cite{huneke}. 
 If 
 $G(I)$ is not  Cohen-Macaulay for an $\m$-primary ideal $I$, then
 ${\rm depth}\,R(I)  ={\rm  depth}\,G(I)+1$ \cite{huc-mar}.  Hence,  
 ${\rm depth}\,R(I)={\rm depth}\,G(I)+1$
 for all $\m$-primary ideals in a regular local  ring.  \qed

A natural question to ask: If the Rees ring is Cohen-Macaulay, then is the fiber cone also. We have a positive answer for $\m$-primary ideals in a two dimensional regular local ring.

\bco
\label{rees-fiber}
Let $(R, \m)$ be a two-dimensional regular local ring. Let $I$ be an $\m$-primary ideal of $R$. If $R(I)$ is Cohen-Macaulay, then $G(I)$ and 
$F(I)$are Cohen-Macaulay.
\eco
\pf The Cohen-Macaulayness of $G(I)$ follows from the Goto-Shimoda theorem. 
In a two dimensional  Cohen-Macaulay local ring, if the Rees ring $R(I)$ is
Cohen-Macaulay, then the reduction number of the ideal $I$ is one by
the Goto-Shimoda theorem. This in turn  implies that  the fiber cone
$F(I)$ is Cohen-Macaulay \cite{shah}.  
\qed

\brm
\been
\item
The converse of Corollary~\ref{rees-fiber} is not true for any arbitrary $\m$-primary ideal. See Example~\ref{ex-2}.

\item
Examples where the associated graded ring is Cohen-Macaulay but the fiber cone is not are well known (\cite[Example~2.8]{heinzer})  and vice versa. 
\eeen
\erm

In this section we give a partial answer to the converse of Corollary~\ref{rees-fiber}.

\bp
\label{rees-fiber-cm}
Let $(R, \m)$ be a two dimensional regular local ring and let $I$ be an  $\m$-primary contracted ideal of $R$. Then 
%\been
%\item
 $F(I)$ is Cohen-Macaulay if and only if $G(I)$ is Cohen-Macaulay.
% 
% \item If $\depth(G(I)) = 1$, then $\depth(F(I)) \geq 1$.
% \eeen
\ep
\pf Since $I$ is contracted we know that if $F(I)$ is Cohen-Macaulay, if and only if  $r(I) \leq 1$ \cite{drv}. Also since $I$ is an ideal in a regular local ring, $G(I)$ is Cohen-Macaulay if and only if $r(I) = 1$. 
Hence $F(I)$ is Cohen-Macaulay if and only of $G(I)$ is Cohen-Macaulay. 

%Now suppose $\depth(G(I)) =1$. Let $a \in I \not \in I^2$ be a superficial element for $I$ whose image in $I/I^2$ is a regular element in $G(I)$. 
%Then  $I^{n+1} :a = I^n$ for all $n \geq 0$. Then from \cite[Proposition~2.7]{cort-zar}, to show that $\depth(F(I)) > 0$,  it is enough to show that $\m I^{n+1} :a = mI^{n}$ for all $n \geq 1$. 

%Put $r = o(I)$. Suppose $c \in \m I^{n+1} :a$. Then $ca \in 
%\m I^{n+1}$. Hence 
%$ca/x^{r(n+1)+1} \in \m I^{n+1}/ x^{r(n+1)+1}$. Hence 
%$c/x^{rn+1} \in {I^{n+1}}^T : a/x^r= {I^{n}}^T$ for all first quadratic transforms $T$ of $I$. This gives $c \in x^{rn + 1}({I^n})^T \cap R = \m I^n$. 
%Hence $\depth(F(I))>0$. If  $\depth(F(I))=2$, then $F(I)$ is Cohen-Macaulay and hence $G(I)$ is Cohen-Macaulay which contradicts  our assumption that $\depth(G(I)) =1$. Hence $\depth(F(I)) =1$.
\qed

To prove the  equality of depth of fiber cone and associated graded ring for contracted ideals, we need the concept of minimal mixed multiplicity.

\section{Ideals of minimal mixed multiplicity  and The Complex $ C(x, \a_{d-1}, {\F}, (1, n))$} 

In this section we prove that for contracted ideals in a two dimensional regular local ring, the fiber cone and the associated graded ring have the same depth. We first show that in a two dimensional regular local ring, complete ideals are precisely 
examples of ideals of minimal mixed multiplicity. 
We first show that contracted ideals in a two-dimensional regular local ring are ideals of minimal mixed multiplicity. Next we consider the depth of the associated  graded ring and fiber cone of an ideal of minimal mixed multiplicity in and Cohen-Macaulay ring of dimension at least two. 

We show that if $I$ is an ideal of  minimal mixed multiplicity in a Cohen-Macaulay local ring of dimension at least two and if the associated graded ring is Cohen-Macaulay then the fiber cone has almost maximal depth. 
We  exploit the graded complex $ C(x, \a_{d-1}, {\F}, (1, n))$, where $x \in \m$ and  $\a_{d-1}= a_1, \ldots, a_{d-1} \in I$.

Let $I$ be an $\m$-primary ideal in a local ring of dimension $d$. It was proved in \cite{bhat}, that for large values of $r$ and $s$, 
$\ell(R/
{\m}^r I^s)$ can be written as a polynomial in $r$ and $s$ and can be written in the form 
\beqn
\sum_{i+j \leq d} {e_{ij}} {r+ i \choose i}{s + j \choose j }
\eeqn
where $e_{ij}$ are integers. When $i+j= d$, the integers we put $e_{ij}=e_j(\m|I)$ for $j = 0, \ldots,d$. In this case these are called the mixed multiplicities of $\m$ and $I$.

\bdef
\cite[Definition~2.2]{drv}
We say that an ideal has minimal mixed multiplicity if $\mu(I) = e_{d-1}(\m|I) + d-1$.
\edeff

 \bdef
  A
set of elements $x_1, \ldots, x_d$ is a called a {\em joint reduction} of
a set of ideals $I_1, \ldots, I_d$ if $x_i \in I_i$ for $i=1, \ldots, d$
and there exists a positive integer $n$ so that \beqn
  \left[ \sum_{j=1}^{d} x_j \, I_1 \cdots \hat{I_j} \cdots I_d \right]
   (I_1 \cdots I_d)^{n-1}
=  (I_1 \cdots I_d)^n.
\eeqn
\edeff

Rees proved that if $R/m$ is infinite, then joint reductions
exist.
 From a result of Rees, it follows that  $e_j(\m|I)$ is the multiplicity of any joint reduction of  
the multi-set of ideals consisting of $j$ copies of $I$ and $d-j$ copies of
$m.$
We shall denote such a multi-set  of ideals by $(\m^{[d-j]}|I^{[j]})$.

\bl
\label{cont-mmm}
Let $I$ be an $\m$-primary ideal in a two dimensional regular local ring $(R, \m)$ . Then 
$I$ is a contracted ideal if and only of $I$ is an ideal if 
 of minimal mixed multiplicity. 
\el
\pf
Every contracted ideal in a regular local ring $(R, \m)$ 
satisfies the property that $\mu(I) = 1 + o(I)$, where $o(I)$ is the $\m$-adic order of $I$ and $o(I) = e(\m| I)$  \cite[Theorem~4.1]{verma1}. 

Conversely, suppose $I$ is an ideal of minimal mixed multiplicity, then there exists a joint reduction $(x, a)$ of $(\m , I)$ such that $\m I = x I + a \m$. By \cite{huneke}, to show that $I$ is contracted it is enough to show that  $\m I : (x) = I$.   

Clearly, $I \subseteq \m I : (x)$. Now suppose $\alpha \in \m I : (x)$. Then
$
\alpha x \in \m I = x I + a \m.
$
Hence $
\alpha x = x i + a \beta
$ 
for some $i \in I$ and $\beta \in \m.
$
Therefore
$
x(\alpha - i) = a \beta.
$
Since $(x,a)$ is a regular sequence in $R$, 
$
\alpha - i = a t
$
for some $t \in R$. Hence
$
\alpha = i + at \in I
$ 
since $i, a \in I.
$
\qed

We first prove some results on ideals of minimal mixed multiplicity. 

\bt
 \label{fiber-cm}
 Let $I$ be an ideal of minimal mixed multiplicity in Cohen-Macaulay local ring of dimension $d \geq 2$. If $F(I)$ is Cohen-Macaulay, then $G(I)$ are $R(I)$ are Cohen-Macaulay.
 \et
 \pf From \cite[Corollary~2.5]{drv}, $F(I)$ is Cohen-Macaulay if and only if $r(I) \leq 1$. If $r(I) \leq 1$, then both $G(I)$ is Cohen-Macaulay and hence $R(I)$ is Cohen-Macaulay. 
 \qed

The converse of Theorem~\ref{fiber-cm} need not always hold true, (see Proposition~\ref{depth}, Example~\ref{ex-1}).
 To prove this  we consider a complex 
$C( \a_{d-1}, {\F}, (1, n))$
which will link the depth of the associated graded ring with the depth of the fiber cone. 

Let ${\F} = \{\m  I^{n} \}_{n \geq 0}$  be the filtration and let $x \in \m$ and $\a_{d-1} \in I$. The complex  
$C(\a_{d-1}, {\F}, (0, n))$,
 $C( \a_{d-1}, {\F}, (1, n))$
$C(x, \a_{d-1}, {\F}, (1, n))$ can be constructed via the mapping cylinder. For 
details on these complexes 
we request the reader to refer to \cite {huc-mar} and \cite{anna-clare-2}.
 We have the short exact sequence:
 
 \mptmcomx
\beq
\label{complex-1}
      0 
\lrar C(\a_{d-1}, {\F}, (1, n))
\lrar C(x, \a_{d-1}, {\F}, (1, n))
\lrar C(\a_{d-1}, {\F}, (0, n))
\lrar 0
\eeq
where 
\beqn
\begin{array}{llllllllll}

       C(\a_{d-1}, {\F}, (0, n))
:    &  0 
\rar& \f{R}{I^{n-d+1}} 
&\rar \cdots 
\rar \left( {\f{R}{I^{n-d+i}}} \right)^{d-1 \choose i} 
&\rar \cdots
&\rar& \f{R}{I^n} 
&\rar& 0\\ 
      C(\a_{d-1}, {\F}, (1, n))
:    &   0 
\rar& \f{R}{\m I^{n-d+1}}  
&\rar \cdots 
\rar \left(
            \f{R}
              {\m I^{n-d+i}} 
             \right)^{d-1 \choose i} 
&\rar  \cdots 
&\rar &\f{R}{\m I^n} 
&\rar& 0\\
     % \end{array}
      %\eeqn    
%\beqn
     C(x, \a_{d-1}, {\F}, (1, n)) 
 :   &   0 
\rar& \f{R}{I^{n-d+1}}     
&\rar \cdots  
\rar  { \displaystyle \bigoplus_{j=0}^{1} }
     \left( 
     \f{R}
       {\m^{1-j} I^{n-d + i-j}} 
     \right)^{{d-1 \choose i-j}}   
     &\rar 
     \cdots       
&\rar& \f{R}{\m I^{n}}  
&\rar&    0.\\
\end{array}
 \eeqn

Let $H_i(-)$  denote the $i$-th homology of the corresponding complex.
 From \ref{complex-1} we have the long exact sequence of homologies
 \beq
\label{complex-2}
\begin{array}{lclclc}
       \cdots
&\rar& H_{i+1} (C(x, \a_{d-1}, {\F}, (1, n)))
&\rar& H_{i} (C(\a_{d-1}, {\F}, (0, n)))
&\rar \\
       H_{i} (C(\a_{d-1}, {\F}))(1, n)
&\rar& H_{i} (C(x, \a_{d-1}, {\F}, (1, n)))
&\rar& H_{i-1} (C(\a_{d-1}, {\F}, (0, n)))
& \cdots .
\end{array}
\eeq

We also have  exact sequence:
 \beqn
         0 
\lrar K.(\a_{d-1}^{o},  F(I))( n)
\lrar C( \a_{d-1},   {\F},  (1, n))
\lrar C(\a_{d-1}, {\F},  (0, n))
\lrar 0
\eeqn
gives rise to the long exact sequence:
\beq
\label{vanishing-d-1-fiber}
\begin{array}{lclclc}
       \cdots
&\rar& H_{i+1} (C(\a_{d-1}, {\F}, (1, n)))
&\rar& H_{i+1} (C(\a_{d-1}, {\F}, (0, n)))
&\rar \\
       H_{i} (K. (\a_{d-1},  F(I))( n)
&\rar& H_{i} (C(\a_{d-1}, {\F}, (1, n)))
&\rar& H_{i} (C(\a_{d-1}, {\F}, (0, n)))
&  .
\end{array}
\eeq

\mptm
 
 \bl
\label{lemma-vanishing}
Let $I$ be an $\m$-primary ideal of minimal mixed multiplicity in a Cohen-Macaulay local ring $(R, \m)$.
 Let $x \in \m$ and 
 $\a_{d-1} \in I$ 
 be  a joint reduction of $\m$ and $I$. Then
 $H_i( C(x, \a_{d-1},  {\F}, (1, n)) ) = 0$ for all $n \geq 1$ and for all $i \geq 1$. 
 \el
 \pf It is enough to show that  $H_1( C(x, \a_{d-1},  {\F}, (1, n)) ) = 0$ for all $n \geq 1$ \cite{anna-clare-2}. From \cite{anna-clare-2},
 \beqn
 H_1( C(x, \a_{d-1},  {\F}, (1, n)) ) 
 \cong \f{(x, \a_{d-1}) \cap \m I^n}{x I^n + \a_{d-1} \m I^{n-1}}, 
 \hspace{.5in} \mbox{for all } n \geq 1.
  \eeqn
 Since $I$ is an ideal of minimal mixed multiplicity, for all $n \geq 1$, 
 $
 \m I^n= x I^n + \a_{d-1}^n \m
$. Hence 
  \beqn
        (x, \a_{d-1}) \cap \m I^n
 \cong (x, \a_{d-1})  \cap (ex I^n + \a_{d-1}^n \m)
 =      x I^n + \a_{d-1}^n \m
 =      \m I^n.
 \eeqn
 Therefore
 $
 H_1( C(x, \a_{d-1},  {\F}, (1, n)) ) = 0$,
  for all $n \geq 1$.
 \qed
 
% \bp
% \label{prop-vanishing}
% Let $I$ be an $\m$-primary ideal of minimal mixed multiplicity in a Cohen-Macaulay ring of dimension $d \geq 2$. Let
% $x, \a_{d-1}$ be a joint reduction of  
% $(\m, I^{d-1})$ such that 
% $\m I  = x I + \a_{d-1} \m$. 
% \beq
%\label{vanishing-d-1}
%H_{i} (C(\a_{d-1},      {\F}, (0, n)))
%\cong
%       H_{i} (C(\a_{d-1},      {\F}, (1, n))),
%       \hspace{.5in} \mbox{ for all } i \geq 1, n \geq 0.
%       \eeq
%\ep
% \pf From the complex (\ref{complex-1})
% we have the long exact sequence of homologies:
%\beq
%\begin{array}{lclclc} \nno
%         \cdots
%&\rar& H_{i+1}   (C(x, \a_{d-1}, {\F}, (1, n)))
%&\rar& H_{i} (C(\a_{d-1},      {\F}, (0, n)))
%&\rar \\
%\label{complex-2}
%       H_{i} (C(\a_{d-1},      {\F}, (1, n)))
%&\rar& H_{i} (C(x, \a_{d-1}, {\F}, (1, n)))
%&\rar& \cdots
%&  .
%\end{array}
%\eeq
%From  Lemma~\ref{lemma-vanishing} and the exact sequence  (\ref{complex-2}) we get
%Hence 
%\beqn
%H_{i} (C(\a_{d-1},      {\F}, (0, n)))
%\cong
%       H_{i} (C(\a_{d-1},      {\F}, (1, n))),
%       \hspace{.5in} \mbox{ for all } i \geq 1, n \geq 0.
%       \eeqn
%\qed

 \bt
 \label{depth}
 Let $I$ be an $\m$-primary ideal of minimal mixed multiplicity in a Cohen-Macaulay ring of dimension $d \geq 2$. 
 Let
 $x, a_{d-1}$ be a joint reduction of  
 $\m, I^{d-1}$. 
 If   $a_{1}^{*},\ldots, a_{t}^{*} $ is a regular sequence in $G(I)$, then  
 then $a_{1}^{o},\ldots, a_{t}^{o} $
 $\depth(F(I)) \geq d-1$. In particular, 
If $\depth(G(I)) \geq d-1 $ , then $\depth(F(I)) \geq d-1$.
  \et
 \pf  We can choose $x \in \m$ and $a_1, \ldots, a_{d-1} \in I$ to be a superficial sequence for $\m$ and $I$, 
 so that $a_{d-1}^{*}$ is a superficial sequence in $G(I)$ and $F(I)$. Since 
 $a_{1}^{*},\ldots, a_{t}^{*} $ is a regular sequence in 
 $G(I)$ we have
 %%%%%%%%%%%%%%%%%%%%%
 from  \cite[Proposition~3.3]{huc-mar}, we get
 \beq
 \label{van-1}
   H_i \left( C(\a_{t}, {\F}, (0, n)) \right) = 0, 
   \hspace{.5in} & \mbox{ for all } 
 i \geq d-t+1 \hspace{.5in} &\cite[Proposition~3.3]{huc-mar}\\
 \label{van-2}
  H_{i} (C(\a_{t},      {\F}, (1, n))),
       \hspace{.5in} & \mbox{ for all } i \geq d-t+ ,  n \geq 0
       & \mbox{from (\ref{complex-2})} \\
       H_{i} ( K.(\a_{t}, F(I))( n) = 0, \hspace{.5in} 
       &\mbox{ for all } i \geq d-t+1, n \geq 0
       & \mbox{applying  (\ref{van-1})  and (\ref{van-2}) 
to (\ref{vanishing-d-1-fiber})}.
 \eeq
\qed

\bt
\label{depth-fiber}
Let $I$ be an ideal of minimal mixed multiplicity in a Cohen-Macaulay local ring of dimension two. Let  $(x, a)$ be a joint reduction of $I$. 
If $a^{o}$ is a regular element in $F(I)$, then $a^{*}$ is a regular element in $G(I)$. In particular if $\depth(F(I)) \geq 1$, then 
$\depth(G((I)) \geq 1$. 
\et
\pf Since $I$ is an ideal of minimal mixed multiplicity, from Lemma~\ref{lemma-vanishing},    $H_i( C(x, a,  {\F}, (1, n)) ) = 0$ for all $n \geq 1$ and for all $i \geq 1$.  Hence from (\ref{complex-2})
\beq
\label{vanishing-d-1-1}
H_{1} (C(a,      {\F}, (0, n)))
\cong
       H_{1} (C(a,      {\F}, (1, n))),
       \hspace{.5in} \mbox{ for all }n \geq 0.
       \eeq
Since  $a^{o}$ is a regular element in $F(I)$, from (\ref{vanishing-d-1-fiber})
\beq
\label{vanishing-d-1}
0 &=&  H_{0} (C(a^{o},  F(I),  n))\\ \nno
 H_{1} (C(a,      {\F}, (0, n)))
&\cong&
       H_{1} (C(a,      {\F}, (1, n))),
       \hspace{.5in} \mbox{ for all }n \geq 0.
       \eeq
  and that  the exact sequence:
\beqn
     0 
\rar H_{0} (C(a^{o},  F(I),  n))
\rar H_{0} (C(a, {\F}, (1, n)))
\rar H_{0} (C(a, {\F}, (0, n)))
\rar 0,
\eeqn
is exact for all $n \geq 0$, i.e.
\beqn
      0
\rar \f{I^n}{aI^{n-1} + \m I^n}
\rar \f{R}{a + \m I^n}
\rar \f{R}{a + I^n}
\rar 0, \hspace{.2in} \mbox{ for all }  n \geq 0.
\eeqn
Therefore 
\beqn
         \f{I^n}{aI^{n-1} + \m I^n}
\cong \f{a + I^n}{a + \m I^n}
\cong \f{I^n}{(a \cap I^n) + \m I^n}. 
\eeqn
This implies that 
\beqn
           (a \cap I^n) 
\subseteq a I^{n-1} + \m I^n
=         a I^{n-1}  + x I^{n}. 
\eeqn
The last inequality follows since $I$ is an ideal of minimal mixed multiplicity.
We claim that: 
\beqn
(a \cap I^n) = a I^{n-1} + a x (I^n :a).
\eeqn
Suppose $\alpha \in (a \cap I^n)$. Then $\alpha = a \beta = a i + xj$, where 
$\beta \in R$, $i \in I^{n-1}$ and $j \in I^n$. 
This implies that 
\beqn
a(\beta - i) = xj.
\eeqn
Since $x,a$ is a regular sequence in $R$, $j = aj^{\prime}$ where
$j^{\prime} \in (I^n :a)$.
Therefore 
\beqn
(a \cap I^n) \subseteq a I^{n-1} + ax (I^n :a).
\eeqn
Clearly $a I^{n-1} + ax (I^n :a) \subseteq (a \cap I^n)$. Hence equality holds. This proves the claim.

Since $(a \cap I^n) = a (I^n :a)$ Hence we have
\beqn
a (I^n:a)  \subseteq a I^{n-1} + ax (I^n :a).
\eeqn
and   by Nakayama's lemma, 
\beqn
a (I^n:a)  = a I^{n-1}.
\eeqn
As $a$ is a regular element in $R$, 
\beqn
 (I^n:a)  =  I^{n-1}, \hspace{.5in} 
\mbox{ for all } n \geq 0.
\eeqn
Therefore $a^{*}$ is a regular element in $G(I)$. Hence $\depth ~G(I) \geq 1$. 
\qed

\bt
\label{rees-graded-fiber}
Let $I$ be a contracted ideal  in a two dimensional regular local ring. Then $\depth(R(I)) -1 = \depth(G(I)) =\depth(F(I))$.
\et
 \pf  
 By Theorem~\ref{rees-graded}, $\depth(R(I)) = \depth(G(I))-1$. Hence we only need to prove that $\depth(G(I)) =\depth(F(I))$. By Corollary~\ref{rees-fiber-cm}, $G(I)$ is Cohen-Macaulay if and only $F(I)$ is Cohen-Macaulay. 
 
% Now suppose $\depth(G(I))=1$, then 
% is Cohen-Macaulay, then by \cite[Lemma~1.2]{rees} we can choose a joint reduction $(x,a)$ so that $a^{*}$ is a nonzerodivisor in $G(I)$. Then by Corollary~\ref{rees-fiber-cm}, $\depth(F(I)) \geq 1$ is Cohen-Macaulay. 
% Clearly, $F(I)$ cannot be Cohen-Macaulay, because, by Corollary~\ref{rees-fiber-cm} it would imply that $G(I)$ is. Similarly, if $\depth(F(I)) = 1$, then by
%Theorem~\ref{depth-fiber}, and Corollary~\ref{rees-fiber-cm} we can conclude $\depth(F(I)) = 1$.

By Theorem~\ref{rees-fiber-cm} and Theorem~\ref{depth-fiber} we get
	$\depth(G(I)) = 1$ if and only if $\depth(F(I)) = 1$.
	The only case remaining is when $\depth(G(I)) = 0$. In this case, $\depth(F(I))=0$ and vice versa.
\qed

\section{Some Examples}
We present an infinite class of examples where the Rees ring is Cohen-Macaulay but the fiber cone is not.

%\bex
%\label{ex-1}
%{\em 
%Let $R = k[X, Y, Z]_{\m}$, where $\m = (ex,Y, Z)$. Let\\
%$I = (X^4, X^3Y, X^2Z, X^2Y^2, XY^2Z, XYZ^2, XZ^3, Y^3, Y^2Z^2, YZ^3,
%Z^5)$. }
%\eex
%Then $J= (X^4 + YZ^3, X^2Z, Y^3 + Z^5 )$ is a minimal reduction of
%$I$. The following hold true:\\
%(i) $J \cap I^2 = JI$\\
%(ii) $JI^n = I^{n+1}$ for all $n \geq 2 $. 

%Hence $R(I)$ and $G(I)$ is Cohen-Macaulay \cite{goto-s},
%\cite{huc-mar}. 

%It is easy to see that 
%$
%\mu(I) := \ell(I/ \m I) = 11 = 9+2 =e(X+Z, X^2Z, Y^3 + Z^5 ) + 2.
%$ 
%Hence $I$ is an ideal of minimal mixed multiplicity and 
%the Hilbert  series of $F(I)$ is of the from \cite[Theorem~2.4]{drv}
%$$
%H(F(I), t) = \f{1+8t}{(1-t)^3}
%$$
%Since $I$ is an ideal of minimal mixed multiplicity and the reduction
%number of $I$, $r(I) >1$, 
%$F(I)$ is not Cohen-Macaulay \cite[Corollary~2.5]{drv}. 
%By Theorem~\ref{depth},  depth~$F(I) =2$.
%\qed

%We now present an infinite class of examples where the fiber cone is
%Cohen-Macaulay but the Rees ring is not. 
\bex
\label{ex-1}
{\em 
Let $R = k[X, Y, Z]_{\m}$, where $\m = (X,Y, Z)$. Let $s \geq 2$ and \\
$I_s = (
 Y^3, 
     X^{s+1}Y^2, 
     Y^2XZ, 
     Y^2Z^{3}, 
     X^{2s+1}Y,
     YZ^2X,   
      YZ^{5}, 
      X^{3s+1}, 
     ZX^2,  
     Z^4X,
     Z^{7})
$. 

}
\eex
Then $J= (X^{3s+1} + Z^7, X^2Z, Y^3 )$ is a minimal reduction of
$I$. The following hold true:\\
(i) $J \cap I^2 = JI$\\
(ii) $JI^n = I^{n+1}$ for all $n \geq 2 $. 

Hence $R(I)$ and $G(I)$ is Cohen-Macaulay \cite{goto-s},
\cite{huc-mar}. 

It is easy to see that 
$
\mu(I) := \ell(I/ \m I) = 11 = 9+2 =e(X+Y+Z, X^2Z, Y^3 + Z^7) + 2.
$ 
Hence $I$ is an ideal of minimal mixed multiplicity and 
the Hilbert  series of $F(I)$ is of the from \cite[Theorem~2.4]{drv}
$$
H(F(I), t) = \f{1+8t}{(1-t)^3}
$$
Since $I$ is an ideal of minimal mixed multiplicity and the reduction
number of $I$, $r(I) >1$, 
$F(I)$ is not Cohen-Macaulay \cite[Corollary~2.5]{drv}. 
By Proposition~\ref{depth},  depth~$F(I) =2$.
\qed

We now present an infinite class of examples where the fiber cone is
Cohen-Macaulay but the Rees ring is not. 

\bex
\label{ex-2}
{\em 
Let $R = k[X, Y]_{\m}$, where $\m = (X,Y)$. Let
$$
I = (X^n, X^{2}Y^{n-2},  XY^{n-1}, Y^n).
$$
Then $F(I)=k[X^n, X^{2}Y^{n-2},  XY^{n-1}, Y^n] $  
is Cohen-Macaulay but $R(I)$ and $G(I)$ are not. Moreover if $n$ is
even then $F(I)$ is Gorenstein.   }
\eex

One can verify that
\beqn
            F(I) 
&\cong&  \left\{ \begin{array}{ll}
         {\displaystyle  \f{k[a,b,c,d]}
                           {(c^2 - bc, b^{r+1}- acd^{r-1}, b^r c - ad^r)}}
     & \mbox{ if }  n= 2r+1 \\
      {\displaystyle   
        \f{k[a,b,c,d]}{(c^2 - bc,  b^{r}- ad^{r-1})}}     
& \mbox{ if }  n = 2r \\       
                                           \end{array}
                                           \right. .
%}
\eeqn        

Clearly, for all $n \geq 4$, the image of $a,d$ in $F(I)$ is a regular
sequence. Hence $F(I)$ is Cohen-Macaulay. 

Moreover,  the Hilbert-Series of $F(I)$ is 

 \beqn  
H(F(I), t)  &=& \left\{ \begin{array}{ll}
            {\displaystyle  \f{1 + 2(t + \cdots + t^r)}{(1-t)^2}}
     & \mbox{ if }  n= 2r+1 \\
      {\displaystyle   \f{1 + 2(t + \cdots + t^{r-1}) +
          t^{r}}{(1-t)^2}}      & \mbox{ if }  n = 2r \\       
                                           \end{array}
                                           \right. .
%}
\eeqn        

Put $J = (X^n, Y^n)$. 

Put $r = \left\lfloor n/2 \right\rfloor$. 
Then $X^{n-1}  Y^{n (r-1) + 1} \in  I^r \sms J I^{r-1}$. 
and 
 $JI^{r}  = I^{r+1}$ which implies that the reduction number of $I$
 with respect to $J$ is at least two. 
 Now $J \cap I^2 \not = JI$ since  $X^3Y^{2n-3}
 \in J \cap I^2 \sms JI$. Therefore 
 $G(I)$ is not Cohen-Macaulay. 
Also $R(I)$ are not Cohen-Macaulay. 
\qed

%\eex
%Let $I = (X^7, X^4y-XY^4,  X^2Y^6,  Y^8)$, 
% $J = (X^4Y-XY^4, X^7 + Y^8)$. The one can verify that $J \cap I^2 \not = JI$
% but $J \cap \m I^n = \m J I^{n-1}$ for all $n \geq 1$. Hence $F(I)$ is Cohen-Macaulay but $G(I)$ is not Cohen-Macaulay. 

%

%Let $I = (Z^7, X^4y-XY^4,  X^2y^6,  y^8)$, 

%An example of an  ideal in a two dimensional regular local ring
%where depth of $F(I)$ is one but  is  depth of $G(I))$ is zero.

%\bex
%\label{ex-3}
%{\em 
%Let $R= k[X,Y]_{\m}$, where $\m = (X,Y)$. Let $I = (X^4, X^3Y^4, XY^{12}, Y^{16})$. Then $J = (X^4, Y^{16})$ is a minimal reduction of $I$. One can verify that 
%\been
%\item
%$I^2 = JI + (X^6 Y^8)$, $I^3 = J^I2$;

%\item
%$I^2:J = I + (X^2 Y^8)$;

%\item
%$J \cap I^2 = I^2 = JI + (X^6 Y^8)$;
%\eeen

%Hence $\depth(G(I))=0$;

%One can verify that 
%$ (\m I^n :(X^4) ) \cap I^{n-1} \cap \m I^{n-1}$ for all  $n\geq 1$.
% Hence by \cite[Lemma~5.1]{dcruz-verma}, ${X^4}^{o}$ is an nonzerodivisor in $F(I)$
%$\depth(F(I)) \geq 1$.
%Since the Hilbert Series of $F(I)$ is $(1 + 2t + 2t^2 - t^3)  /(1-t)^2$, $F(I)$ is not Cohen-Macaulay. 
%}
%\eex

\end{document}